\documentclass[10pt]{amsart}
\usepackage[psamsfonts]{amssymb}
\usepackage{amsthm,amscd}
\usepackage{type1cm}
\usepackage{color}

\newtheorem{thm}{Theorem}

\newtheorem{lem}{Lemma}
\newtheorem{cor}{Corollary}
\newtheorem{pro}{Proposition}
\theoremstyle{remark}

\theoremstyle{definition}

\allowdisplaybreaks[4]

\title[Classification of endo-commutative curled algebras]{A classification of endo-commutative curled algebras of dimension 2 over a non-trivial field}

\author[S.-E. Takahasi]{Sin-Ei Takahasi}
\author[K. Shirayanagi]{Kiyoshi Shirayanagi}
\author[M. Tsukada]{Makoto Tsukada}
\address[S.-E. Takahasi]{Laboratory of Mathematics and Games\\ Katsushika 2-371\\ Funabashi\\ Chiba 273-0032\\ Japan}
\address[K. Shirayanagi, M. Tsukada]{Department of Information Science\\ Toho University\\ Miyama 2-2-1\\ Funabashi\\ Chiba 274-8510\\ Japan}
\email{sin\_ei1@yahoo.co.jp}
\email[K. Shirayanagi (Corresponding author)]{kiyoshi.shirayanagi@is.sci.toho-u.ac.jp}
\email[M. Tsukada]{tsukada@is.sci.toho-u.ac.jp}

\makeatletter
\@namedef{subjclassname@2020}{\textup{2020} Mathematics Subject Classification}
\makeatother

\subjclass[2020]{Primary 17A30; Secondary 17D99, 13A99}
\keywords{Nonassociative algebras, Endo-commutative algebras, Commutative algebras, Curled algebras, Straight algebras.}

\begin{document}




\maketitle

\begin{abstract}
  An endo-commutative algebra is a nonassociative algebra in which the square mapping preserves multiplication. In this paper, we give a
  complete classification of endo-commutative curled algebras of dimension 2 over an arbitrary non-trivial field, where a curled algebra
  satisfies the condition that the square of any element is a scalar multiple of that element. We list all multiplication tables of the algebras up to isomorphism.
\end{abstract}

\section{Introduction}\label{sec:intro}
Let $A$ be a nonassociative algebra.  The square mapping $x\mapsto x^2$ from $A$ to itself yields various important concepts of $A$.  In fact, if the square mapping of $A$ is surjective, then $A$ is said to be square-rootable, see \cite{2dim-comm-char2,z3}.  Also, as is well known, if the square mapping of $A$ preserves addition, then $A$ is said to be anti-commutative.  Moreover, if the square mapping of $A$ is the zero mapping, then $A$ is said to be zeropotent.  We refer the reader
to \cite{z1,z2,z3} for the details on zeropotent algebras. 

The subject of this paper is another concept that also naturally arises from the square mapping.  We define $A$ to be {\it{endo-commutative}}, if the square mapping of $A$ preserves multiplication, that is, $x^2y^2=(xy)^2$ holds for all $x, y\in A$.
This terminology comes from the identity $(xx)(yy)=(xy)(xy)$ that depicts the {\it innerly} commutative property.
The concept of endo-commutativity was first introduced in \cite{TST1}, where we gave a complete classification of endo-commutative algebras of dimension 2 over the trivial field $\mathbb F_2$ of two elements.

We separate two-dimensional algebras into two categories: {\it{curled}} and {\it{straight}}.  That is, a two-dimensional algebra is curled if the square of any element $x$ is a scalar multiple of $x$, otherwise it is straight.  Research related to curled algebras can be found in \cite{level2,length1}. 

The aim of this paper is to completely classify endo-commutative curled algebras of dimension 2 over an arbitrary non-trivial field.  The strategy for the classification is based on that of Kobayashi \cite{2dim-comm-char2}.  We can find classifications of associative algebras of dimension 2 over the real and complex number fields
in \cite{onkochishin}. For other studies on two-dimensional algebras, see \cite{moduli,ABR,2dim,variety,classification}. 

Our main theorem states that endo-commutative curled algebras of dimension 2 over an arbitrary non-trivial field are classified into three families $C_0, C_1$ and  $\{C_2(a)\}_{a\in K}$ up to isomorphism, where $K$ is the base field and $C_i$ is a type of the structure matrix, that is, the matrix of structure constants with respect to a linear base.  The details are described in Theorem 1. As an application of the main theorem, we further give a complete classification of
endo-commutative curled algebras of dimension 2 over an arbitrary non-trivial field in each case of zeropotent, commutative, anti-commutative and associative algebras.
The details are described in Corollary 4. Note that a result similar to the main theorem appears in Asrorov-Bekbaev-Rakhimov \cite{ABR}, despite the difference in
canonical representation.

\section{A criterion for isomorphism of two-dimensional algebras}\label{sec:iso-criterion}
Let $K$ be a non-trivial field and $GL_n(K)$ the general linear group consisting of nonsingular $n\times n$ matrices over $K$.  For any $X=\small{\begin{pmatrix}a&b\\c&d\end{pmatrix}}\in GL_2(K)$, define
\[
\widetilde{X}=\begin{pmatrix}a^2&b^2&ab&ab\\c^2&d^2&cd&cd\\ac&bd&ad&bc\\ac&bd&bc&ad\end{pmatrix}.
\]
Then we have the following.

\begin{lem}
The mapping $X\mapsto\widetilde{X}$ is a group-isomorphism from $GL_2(K)$ into $GL_4(K)$ with $|\widetilde X|=|X|^4$.
\end{lem}

\begin{proof}
Straightforward.
\end{proof}

Let $A$ be a 2-dimensional algebra over $K$ with a linear base $\{e, f\}$.  Then $A$ is determined by the multiplication table $\begin{pmatrix}e^2&ef\\fe&f^2\end{pmatrix}$.  We write
\[
\left\{
    \begin{array}{@{\,}lll}
     e^2=a_1e+b_1f\\
     f^2=a_2e+b_2f\\
     ef=a_3e+b_3f\\
     fe=a_4e+b_4f\\
   \end{array}
  \right. 
\]
with $a_i, b_i\in K\, \, (1\le i\le4)$ and the matrix
$\begin{pmatrix}a_1&b_1\\a_2&b_2\\a_3&b_3\\a_4&b_4\end{pmatrix}$ is called the structure matrix of $A$ with respect to the base $\{e, f\}$. 

We hereafter will freely use the same symbol $A$ for the matrix and for the algebra because the algebra $A$ is determined by its structure matrix.  Then we have the following.

\begin{pro}
Let $A$ and $A'$ be 2-dimensional algebras over $K$.  Then $A$ and $A'$ are isomorphic iff there is $X\in GL_2(K)$ such that 
\begin{equation}
A'=\widetilde{X^{-1}}AX.
\end{equation}
\end{pro}

\begin{proof}
Let $A$ and $A'$ be the structure matrices of $A$ on a linear base $\{e, f\}$ and $A'$ on a linear base $\{g, h\}$, respectively.  We write
\[
A=\begin{pmatrix}a_1&b_1\\a_2&b_2\\a_3&b_3\\a_4&b_4\end{pmatrix}\, \, {\rm{and}}\, \, A'=\begin{pmatrix}c_1&d_1\\c_2&d_2\\c_3&d_3\\c_4&d_4\end{pmatrix}.
\]

Suppose that $\Phi : A\to B$ is an isomorphism and let $X=\begin{pmatrix}a&b\\c&d\end{pmatrix}\, \, (a, b, c, d\in K)$ be the matrix associated with $\Phi$, that is, $\small{\begin{pmatrix}\Phi(e)\\\Phi(f)\end{pmatrix}=X\begin{pmatrix}g\\h\end{pmatrix}}$.  Since $\Phi$ is a bijection, it follows that $X\in GL_2(K)$.  Also, we have
\begin{equation}
\begin{pmatrix}\Phi(e)^2\\\Phi(f)^2\\\Phi(e)\Phi(f)\\\Phi(f)\Phi(e)\end{pmatrix}
=\begin{pmatrix}\Phi(e^2)\\\Phi(f^2)\\\Phi(ef)\\\Phi(fe)\end{pmatrix}
=A\begin{pmatrix}\Phi(e)\\\Phi(f)\end{pmatrix}=AX\begin{pmatrix}g\\h\end{pmatrix}
\end{equation}
On the other hand, we have
\[
\left\{
    \begin{array}{@{\,}lll}
    \Phi(e)^2=(ag+bh)^2=a^2g^2+b^2h^2+ab(gh+hg)\\
    \Phi(f)^2=(cg+dh)^2=c^2g^2+d^2h^2+cd(gh+hg)\\
     \Phi(e)\Phi(f)=\{ag+bh\}\{cg+dh\} =acg^2+bdh^2+adgh+bchg\\
     \Phi(f)\Phi(e)=\{cg+dh\}\{ag+bh\} =cag^2+dbh^2+cbgh+dahg,
   \end{array}
  \right. 
\]
hence
\begin{equation}
\begin{pmatrix}\Phi(e)^2\\\Phi(f)^2\\\Phi(e)\Phi(f)\\\Phi(f)\Phi(e)\end{pmatrix}=\begin{pmatrix}a^2&b^2&ab&ab\\c^2&d^2&cd&cd\\ac&bd&ad&bc\\ac&bd&bc&ad\end{pmatrix}\begin{pmatrix}g^2\\h^2\\gh\\hg\end{pmatrix}
\end{equation}
By (3) and (2), we have
\begin{equation}
\begin{pmatrix}a^2&b^2&ab&ab\\c^2&d^2&cd&cd\\ac&bd&ad&bc\\ac&bd&bc&ad\end{pmatrix}A'\begin{pmatrix}g\\h\end{pmatrix}=AX\begin{pmatrix}g\\h\end{pmatrix}
\end{equation}
because $\begin{pmatrix}g^2\\h^2\\gh\\hg\end{pmatrix}=A'\begin{pmatrix}g\\h\end{pmatrix}$.  By (4), we have
\[
\widetilde{X}A'=\begin{pmatrix}a^2&b^2&ab&ab\\c^2&d^2&cd&cd\\ac&bd&ad&bc\\ac&bd&bc&ad\end{pmatrix}A'=AX,
\]
and hence we get (1) because $\widetilde{X}^{-1}=\widetilde{X^{-1}}$ from Lemma 1.  

Conversely, suppose that there is $X\in GL_2(K)$ satisfying (1).  Let $\Phi : A\to A'$ be the linear mapping defined by $\small{\begin{pmatrix}\Phi(e)\\\Phi(f)\end{pmatrix}=X\begin{pmatrix}g\\h\end{pmatrix}}$.  Then we can easily see that $\Phi$ is isomorphic by following the reverse of the above argument, hence $A$ and $A'$ are isomorphic.
\end{proof}
\vspace{2mm}

\begin{cor} Let $A$ and $A'$ be two-dimensional algebras over $K$.  If $A$ and $A'$ are isomorphic, then ${\rm{rank}}\, A={\rm{rank}}\, A'$.
\end{cor}
\vspace{2mm}

When (1) holds, we say that the matrices $A$ and $A'$ are equivalent and refer $X$ as a $\it{transformation\, \, matrix}$ for the equivalence $A\cong A'$.  Also, we call this $X$ a transformation matrix for the isomorphism $A\cong A'$ as well.
\vspace{5mm}

\section{Two-dimensional endo-commutative algebras}\label{sec:ec-algebras}
Let $A$ be a 2-dimensional endo-commutative algebra over a non-trivial field $K$ with structure matrix $A=\begin{pmatrix}a_1&b_1\\a_2&b_2\\a_3&b_3\\a_4&b_4\end{pmatrix}$.  Let $x, y\in A$ be arbitrary, and write $\small{\left\{
    \begin{array}{@{\,}lll}
     x=x_1e+x_2f\\
     y=y_1e+y_2f,
   \end{array}
  \right. }$
and put
\[
\left\{\begin{array}{@{\,}lll}
     A=x_1^2a_1+x_2^2a_2+x_1x_2(a_3+a_4)\\
     B=x_1^2b_1+x_2^2b_2+x_1x_2(b_3+b_4)\\
     C=y_1^2a_1+y_2^2a_2+y_1y_2(a_3+a_4)\\
     D=y_1^2b_1+y_2^2b_2+y_1y_2(b_3+b_4)\\
     E=x_1y_1a_1+x_2y_2a_2+x_1y_2a_3+x_2y_1a_4\\
     F=x_1y_1b_1+x_2y_2b_2+x_1y_2b_3+x_2y_1b_4,
   \end{array}\right. 
\]
hence $x^2=Ae+Bf, y^2=Ce+Df$ and $xy=Ee+Ff$.  Then 
\[
\left\{\begin{array}{@{\,}lll}
x^2y^2=(ACa_1+BDa_2+ADa_3+BCa_4)e+(ACb_1+BDb_2+ADb_3+BCb_4)f\\
(xy)^2=\{E^2a_1+F^2a_2+EF(a_3+a_4)\}e+\{E^2b_1+F^2b_2+EF(b_3+b_4)\}f.
\end{array}\right. 
\]
Therefore we see that $x^2y^2=(xy)^2$ iff
\begin{equation}
\left\{
    \begin{array}{@{\,}lll}
    ACa_1+BDa_2+ADa_3+BCa_4=E^2a_1+F^2a_2+EF(a_3+a_4)\cdots(5-1)\\
    ACb_1+BDb_2+ADb_3+BCb_4=E^2b_1+F^2b_2+EF(b_3+b_4)\cdots(5-2).
         \end{array}
  \right. 
\end{equation}
Put
\[
\left\{
    \begin{array}{@{\,}lll}
    X_1=x_1^2y_1^2\\
    X_2=x_2^2y_2^2\\
    X_3=x_1^2y_2^2\\
    X_4=x_2^2y_1^2\\
    X_5=x_1x_2y_1y_2\\
    X_6=x_1^2y_1y_2\\
    X_7=x_1x_2y_1^2\\
    X_8=x_1x_2y_2^2\\
    X_9=x_2^2y_1y_2.
         \end{array}
  \right. 
\]

\begin{lem}
The nine polynomials $X_i\, \, (1\le i\le9)$ on $K$ are linearly independent.
\end{lem}

\begin{proof}
Straightforward.
\end{proof}

Note that
\begin{align*}
AC&=a_1^2X_1+a_2^2X_2+a_1a_2X_3+a_1a_2X_4+(a_3+a_4)^2X_5+a_1(a_3+a_4)X_6\\
&\quad\quad\quad\quad\quad\quad\quad\quad\quad\quad\quad+a_1(a_3+a_4)X_7+a_2(a_3+a_4)X_8+a_2(a_3+a_4)X_9,
\end{align*}
\begin{align*}
BD&=b_1^2X_1+b_2^2X_2+b_1b_2X_3+b_1b_2X_4+(b_3+b_4)^2X_5+b_1(b_3+b_4)X_6\\
&\quad\quad\quad\quad\quad\quad\quad\quad\quad\quad\quad+b_1(b_3+b_4)X_7+b_2(b_3+b_4)X_8+b_2(b_3+b_4)X_9,
\end{align*}
\begin{align*}
AD&=a_1b_1X_1+a_2b_2X_2+a_1b_2X_3+a_2b_1X_4+(a_3+a_4)(b_3+b_4)X_5+a_1(b_3+b_4)X_6\\
&\quad\quad\quad\quad\quad\quad\quad\quad\quad\quad\quad+b_1(a_3+a_4)X_7+b_2(a_3+a_4)X_8+a_2(b_3+b_4)X_9,
\end{align*}
\begin{align*}
BC&=a_1b_1X_1+a_2b_2X_2+a_2b_1X_3+a_1b_2X_4+(a_3+a_4)(b_3+b_4)X_5+b_1(a_3+a_4)X_6\\
&\quad\quad\quad\quad\quad\quad\quad\quad\quad\quad\quad+a_1(b_3+b_4)X_7+a_2(b_3+b_4)X_8+b_2(a_3+a_4)X_9,
\end{align*}
\begin{align*}
E^2&=a_1^2X_1+a_2^2X_2+a_3^2X_3+a_4^2X_4+2\{(a_1a_2+a_3a_4)X_5+a_1a_3X_6+a_1a_4X_7+a_2a_3X_8+a_2a_4X_9\}
\end{align*}
\begin{align*}
F^2&=b_1^2X_1+b_2^2X_2+b_3^2X_3+b_4^2X_4+2\{(b_1b_2+b_3b_4)X_5+b_1b_3X_6+b_1b_4X_7+b_2b_3X_8+b_2b_4X_9\}
\end{align*}
and
\begin{align*}
EF&=a_1b_1X_1+a_2b_2X_2+a_3b_3X_3+a_4b_4X_4+(a_1b_2+a_2b_1+a_3b_4+a_4b_3)X_5+(a_1b_3+a_3b_1)X_6\\
&\quad\quad\quad\quad\quad\quad\quad\quad\quad\quad\quad\quad+(a_1b_4+a_4b_1)X_7+(a_2b_3+a_3b_2)X_8+(a_2b_4+a_4b_2)X_9
\end{align*}
Then we have
\begin{align*}
 ACa_1+&BDa_2+ADa_3+BCa_4\\
&=\{a_1^3+a_2b_1^2+a_1a_3b_1+a_1a_4b_1\}X_1+\{a_1a_2^2+a_2b_2^2+a_2a_3b_2+a_2a_4b_2\}X_2\\
&\quad+\{a_1^2a_2+a_2b_1b_2+a_1a_3b_2+a_2a_4b_1\}X_3+\{a_1^2a_2+a_2b_1b_2+a_2a_3b_1+a_1a_4b_2\}X_4\\
&\quad+\{(a_1+b_3+b_4)(a_3+a_4)^2+a_2(b_3+b_4)^2\}X_5\\
&\quad+\{(a_1^2+a_4b_1)(a_3+a_4)+(a_2b_1+a_1a_3)(b_3+b_4)\}X_6\\
&\quad+\{(a_1^2+a_3b_1)(a_3+a_4)+(a_2b_1+a_4a_1)(b_3+b_4)\}X_7\\
&\quad+\{(a_1a_2+a_3b_2)(a_3+a_4)+a_2(a_4+b_2)(b_3+b_4)\}X_8\\
&\quad+\{(a_1a_2+a_4b_2)(a_3+a_4)+a_2(a_3+b_2)(b_3+b_4)\}X_9
\end{align*}
and
\begin{align*}
&E^2a_1+F^2a_2+EF(a_3+a_4)\\
&=\{a_1^3+a_2b_1^2+a_1b_1(a_3+a_4)\}X_1+\{a_1a_2^2+a_2b_2^2+a_2b_2(a_3+a_4)\}X_2\\
&\quad+\{a_1a_3^2+a_2b_3^2+a_3b_3(a_3+a_4)\}X_3+\{a_1a_4^2+a_2b_4^2+a_4b_4(a_3+a_4)\}X_4\\
&\quad+\{2a_1(a_1a_2+a_3a_4)+2a_2(b_1b_2+b_3b_4)+(a_3+a_4)(a_1b_2+a_2b_1+a_3b_4+a_4b_3)\}X_5\\
&\quad+\{2a_1^2a_3+2a_2b_1b_3+(a_3+a_4)(a_1b_3+a_3b_1)\}X_6\\
&\quad+\{2a_1^2a_4+2a_2b_1b_4+(a_3+a_4)(a_1b_4+a_4b_1)\}X_7\\
&\quad+\{2a_1a_2a_3+2a_2b_2b_3+(a_3+a_4)(a_2b_3+a_3b_2)\}X_8\\
&\quad+\{2a_1a_2a_4+2a_2b_2b_4+(a_3+a_4)(a_2b_4+a_4b_2)\}X_9.
\end{align*}
Therefore by Lemma 2, we see that the necessary and sufficient condition for (5-1) to  holds for all $x, y\in A$ is
\[
\left\{
    \begin{array}{@{\,}lll}
   a_1^3+a_2b_1^2+a_1a_3b_1+a_1a_4b_1=a_1^3+a_2b_1^2+a_1b_1(a_3+a_4)\\
a_1a_2^2+a_2b_2^2+a_2a_3b_2+a_2a_4b_2=a_1a_2^2+a_2b_2^2+a_2b_2(a_3+a_4)\\
  a_1^2a_2+a_2b_1b_2+a_1a_3b_2+a_2a_4b_1=a_1a_3^2+a_2b_3^2+a_3b_3(a_3+a_4)\\
  a_1^2a_2+a_2b_1b_2+a_2a_3b_1+a_1a_4b_2=a_1a_4^2+a_2b_4^2+a_4b_4(a_3+a_4)\\
  (a_1+b_3+b_4)(a_3+a_4)^2+a_2(b_3+b_4)^2=2a_1(a_1a_2+a_3a_4)+2a_2(b_1b_2+b_3b_4)\\
 \quad\quad\quad\quad\quad\quad\quad\quad\quad\quad\quad\quad\quad\quad\quad\quad\quad\quad\quad\quad\quad\quad+(a_3+a_4)(a_1b_2+a_2b_1+a_3b_4+a_4b_3)\\
 (a_1^2+a_4b_1)(a_3+a_4)+(a_2b_1+a_1a_3)(b_3+b_4)=2a_1^2a_3+2a_2b_1b_3+(a_3+a_4)(a_1b_3+a_3b_1)\\
 (a_1^2+a_3b_1)(a_3+a_4)+(a_2b_1+a_4a_1)(b_3+b_4)=2a_1^2a_4+2a_2b_1b_4+(a_3+a_4)(a_1b_4+a_4b_1)\\
 (a_1a_2+a_3b_2)(a_3+a_4)+a_2(a_4+b_2)(b_3+b_4)=2a_1a_2a_3+2a_2b_2b_3+(a_3+a_4)(a_2b_3+a_3b_2)\\
 (a_1a_2+a_4b_2)(a_3+a_4)+a_2(a_3+b_2)(b_3+b_4)=2a_1a_2a_4+2a_2b_2b_4+(a_3+a_4)(a_2b_4+a_4b_2).
   \end{array}
  \right. 
  \]
Note that the first two equations always hold in the above nine equations.  Similarly, we have
\begin{align*}
 ACb_1+&BDb_2+ADb_3+BCb_4\\
&=\{a_1^2b_1+b_1^2b_2+a_1b_1b_3+a_1b_4b_1\}X_1+\{a_2^2b_1+b_2^3+a_2b_2b_3+a_2b_2b_4\}X_2\\
&\quad+\{a_1a_2b_1+b_1b_2^2+a_1b_2b_3+a_2b_1b_4\}X_3+\{
a_1a_2b_1+b_1b_2^2+a_2b_1b_3+a_1b_2b_4\}X_4\\
&\quad+\{b_1(a_3+a_4)^2+(b_2+a_3+a_4)(b_3+b_4)^2\}X_5\\
&\quad+\{b_1(a_1+b_4)(a_3+a_4)+(b_1b_2+a_1b_3)(b_3+b_4)\}X_6\\
&\quad+\{b_1(a_1+b_3)(a_3+a_4)+(b_1b_2+a_1b_4)(b_3+b_4)\}X_7\\
&\quad+\{(a_2b_1+b_2b_3)(a_3+a_4)+(b_2^2+a_2b_4)(b_3+b_4)\}X_8\\
&\quad+\{(a_2b_1+b_2b_4)(a_3+a_4)+(b_2^2+a_2b_3)(b_3+b_4)\}X_9
\end{align*} 
and
\begin{align*}
&E^2b_1+F^2b_2+EF(b_3+b_4)\\
&=\{a_1^2b_1+b_1^2b_2+a_1b_1(b_3+b_4)\}X_1+\{a_2^2b_1+b_2^3+a_2b_2(b_3+b_4)\}X_2+\{a_3^2b_1+b_2b_3^2+a_3b_3(b_3+b_4)\}X_3\\
&\quad+\{a_4^2b_1+b_2b_4^2+a_4b_4(b_3+b_4)\}X_4\\
&\quad+\{2b_1(a_1a_2+a_3a_4)+2b_2(b_1b_2+b_3b_4)+(b_3+b_4)(a_1b_2+a_2b_1+a_3b_4+a_4b_3)\}X_5\\
&\quad+\{2a_1a_3b_1+2b_1b_2b_3+(b_3+b_4)(a_1b_3+a_3b_1)\}X_6\\
&\quad+\{2a_1a_4b_1+2b_1b_2b_4+(b_3+b_4)(a_1b_4+a_4b_1)\}X_7\\
&\quad+\{2a_2a_3b_1+2b_2^2b_3+(b_3+b_4)(a_2b_3+a_3b_2)\}X_8\\
&\quad+\{2a_2a_4b_1+2b_2^2b_4+(b_3+b_4)(a_2b_4+a_4b_2)\}X_9.
\end{align*}
Therefore by Lemma 2, we see that the necessary and sufficient condition for (5-2) to holds for all $x, y\in A$ is 
\[
\left\{
    \begin{array}{@{\,}lll}  a_1^2b_1+b_1^2b_2+a_1b_1b_3+a_1b_4b_1=a_1^2b_1+b_1^2b_2+a_1b_1(b_3+b_4)\\
a_2^2b_1+b_2^3+a_2b_2b_3+a_2b_2b_4= a_2^2b_1+b_2^3+a_2b_2(b_3+b_4)\\
a_1a_2b_1+b_1b_2^2+a_1b_2b_3+a_2b_1b_4=a_3^2b_1+b_2b_3^2+a_3b_3(b_3+b_4)\\
a_1a_2b_1+b_1b_2^2+a_2b_1b_3+a_1b_2b_4= a_4^2b_1+b_2b_4^2+a_4b_4(b_3+b_4)\\
b_1(a_3+a_4)^2+(b_2+a_3+a_4)(b_3+b_4)^2\\
\quad=2b_1(a_1a_2+a_3a_4)+2b_2(b_1b_2+b_3b_4)+(b_3+b_4)(a_1b_2+a_2b_1+a_3b_4+a_4b_3)\\
b_1(a_1+b_4)(a_3+a_4)+(b_1b_2+a_1b_3)(b_3+b_4)\\
\quad=2a_1a_3b_1+2b_1b_2b_3+(b_3+b_4)(a_1b_3+a_3b_1)\\
b_1(a_1+b_3)(a_3+a_4)+(b_1b_2+a_1b_4)(b_3+b_4)\\
\quad=2a_1a_4b_1+2b_1b_2b_4+(b_3+b_4)(a_1b_4+a_4b_1)\\
(a_2b_1+b_2b_3)(a_3+a_4)+(b_2^2+a_2b_4)(b_3+b_4)\\
\quad=2a_2a_3b_1+2b_2^2b_3+(b_3+b_4)(a_2b_3+a_3b_2)\\
(a_2b_1+b_2b_4)(a_3+a_4)+(b_2^2+a_2b_3)(b_3+b_4)\\
\quad=2a_2a_4b_1+2b_2^2b_4+(b_3+b_4)(a_2b_4+a_4b_2).
    \end{array}
  \right. 
\]
Note that the first two equations always hold in the above nine equations.  Therefore, the necessary and sufficient condition for (5) to holds for all $x, y\in A$ is 
\begin{equation}
\left\{\begin{array}{@{\,}lll}   a_1^2a_2+a_2b_1b_2+a_1a_3b_2+a_2a_4b_1=a_1a_3^2+a_2b_3^2+a_3b_3(a_3+a_4)\\
  a_1^2a_2+a_2b_1b_2+a_2a_3b_1+a_1a_4b_2=a_1a_4^2+a_2b_4^2+a_4b_4(a_3+a_4)\\
  (a_1+b_3+b_4)(a_3+a_4)^2+a_2(b_3+b_4)^2=2a_1(a_1a_2+a_3a_4)+2a_2(b_1b_2+b_3b_4)\\
 \quad\quad\quad\quad\quad\quad\quad\quad\quad\quad\quad\quad\quad\quad\quad\quad\quad\quad\quad\quad\quad\quad+(a_3+a_4)(a_1b_2+a_2b_1+a_3b_4+a_4b_3)\\
 (a_1^2+a_4b_1)(a_3+a_4)+(a_2b_1+a_1a_3)(b_3+b_4)=2a_1^2a_3+2a_2b_1b_3+(a_3+a_4)(a_1b_3+a_3b_1)\\
 (a_1^2+a_3b_1)(a_3+a_4)+(a_2b_1+a_4a_1)(b_3+b_4)=2a_1^2a_4+2a_2b_1b_4+(a_3+a_4)(a_1b_4+a_4b_1)\\
 (a_1a_2+a_3b_2)(a_3+a_4)+a_2(a_4+b_2)(b_3+b_4)=2a_1a_2a_3+2a_2b_2b_3+(a_3+a_4)(a_2b_3+a_3b_2)\\
 (a_1a_2+a_4b_2)(a_3+a_4)+a_2(a_3+b_2)(b_3+b_4)=2a_1a_2a_4+2a_2b_2b_4+(a_3+a_4)(a_2b_4+a_4b_2)\\
a_1a_2b_1+b_1b_2^2+a_1b_2b_3+a_2b_1b_4=a_3^2b_1+b_2b_3^2+a_3b_3(b_3+b_4)\\
a_1a_2b_1+b_1b_2^2+a_2b_1b_3+a_1b_2b_4= a_4^2b_1+b_2b_4^2+a_4b_4(b_3+b_4)\\
b_1(a_3+a_4)^2+(b_2+a_3+a_4)(b_3+b_4)^2\\
\quad=2b_1(a_1a_2+a_3a_4)+2b_2(b_1b_2+b_3b_4)+(b_3+b_4)(a_1b_2+a_2b_1+a_3b_4+a_4b_3)\\
b_1(a_1+b_4)(a_3+a_4)+(b_1b_2+a_1b_3)(b_3+b_4)\\
\quad=2a_1a_3b_1+2b_1b_2b_3+(b_3+b_4)(a_1b_3+a_3b_1)\\
b_1(a_1+b_3)(a_3+a_4)+(b_1b_2+a_1b_4)(b_3+b_4)\\
\quad=2a_1a_4b_1+2b_1b_2b_4+(b_3+b_4)(a_1b_4+a_4b_1)\\
(a_2b_1+b_2b_3)(a_3+a_4)+(b_2^2+a_2b_4)(b_3+b_4)\\
\quad=2a_2a_3b_1+2b_2^2b_3+(b_3+b_4)(a_2b_3+a_3b_2)\\
(a_2b_1+b_2b_4)(a_3+a_4)+(b_2^2+a_2b_3)(b_3+b_4)\\
\quad=2a_2a_4b_1+2b_2^2b_4+(b_3+b_4)(a_2b_4+a_4b_2).    
\end{array} \right. 
\end{equation}
Put
\[
\left\{\begin{array}{@{\,}lll}   
a=a_1\\
b=b_1\\
c=a_2\\
d=b_2\\
e=a_3\\
f=b_3\\
g=a_4\\
h=b_4
\end{array} \right. 
\]
Then (6) can be rewritten as
\begin{equation}
\left\{\begin{array}{@{\,}lll}   
a^2c+bcd+ade+bcg=ae^2+cf^2+e^2f+efg \cdots(\rm{i})\\
  a^2c+bcd+bce+adg=ag^2+ch^2+egh+g^2h\cdots(\rm{ii})\\
  (a+f+h)(e+g)^2+c(f+h)^2=2a(ac+eg)+2c(bd+fh)+(e+g)(ad+cb+eh+gf)\cdots(\rm{iii})\\
 (a^2+gb)(e+g)+(cb+ae)(f+h)=2a^2e+2cbf+(e+g)(af+eb)\cdots(\rm{iv})\\
 (a^2+eb)(e+g)+(cb+ga)(f+h)=2a^2g+2cbh+(e+g)(ah+gb)\cdots(\rm{v})\\
 (ac+ed)(e+g)+c(g+d)(f+h)=2ace+2cdf+(e+g)(cf+ed)\cdots(\rm{vi})\\
 (ac+gd)(e+g)+c(e+d)(f+h)=2acg+2cdh+(e+g)(ch+gd)\cdots(\rm{vii})\\
abc+bd^2+adf+bch=be^2+df^2+ef^2+efh\cdots(\rm{viii})\\
abc+bd^2+bcf+adh= bg^2+dh^2+fgh+gh^2\cdots(\rm{ix})\\
b(e+g)^2+(d+e+g)(f+h)^2=2b(ac+eg)+2d(bd+fh)+(f+h)(ad+cb+eh+gf)\cdots(\rm{x})\\
b(a+h)(e+g)+(bd+af)(f+h)=2aeb+2bdf+(f+h)(af+eb)\cdots(\rm{xi})\\
b(a+f)(e+g)+(bd+ah)(f+h)=2agb+2bdh+(f+h)(ah+gb)\cdots(\rm{xii})\\
(cb+df)(e+g)+(d^2+ch)(f+h)=2ceb+2d^2f+(f+h)(cf+ed)\cdots(\rm{xiii})\\
(cb+dh)(e+g)+(d^2+cf)(f+h)=2cgb+2d^2h+(f+h)(ch+gd)\cdots(\rm{xiv}).    
\end{array} \right. 
\end{equation}
\vspace{3mm}

(I) We can easily see that (i) and (ii) implies (iii).
\vspace{3mm}

(II) We can easily see that (viii) and (ix) imply (x).
\vspace{3mm}

Moreover, note that
 
 (iv) $\Leftrightarrow a^2g+bg^2+bch+aeh=a^2e+bcf +be^2+afg$,

 (v) $\Leftrightarrow a^2e+be^2+bcf+afg=a^2g+bch+aeh+bg^2$,

 (vi) $\Leftrightarrow acg+cgh+cdh=ace+cdf+cef$,

 (vii) $\Leftrightarrow ace+cef+cdf=acg+cdh+cgh$,

 (xi) $\Leftrightarrow abg+bgh+bdh=abe+bdf+bef$,

 (xii) $\Leftrightarrow abe+bef+bdf=abg+bdh+bgh$,

 (xiii) $\Leftrightarrow bcg+dfg+d^2h+ch^2=bce+d^2f+cf^2+deh$,

 (xiv) $\Leftrightarrow bce+deh+d^2f+cf^2=bcg+d^2h+dfg+ch^2$.

\noindent
By (I), (II) and the above notes, (7) can be rewritten as
\begin{equation}
\left\{\begin{array}{@{\,}lll}   
a^2c+bcd+ade+bcg=ae^2+cf^2+e^2f+efg \cdots(\rm{i})\\
  a^2c+bcd+bce+adg=ag^2+ch^2+egh+g^2h\cdots(\rm{ii})\\
a^2g+bg^2+bch+aeh=a^2e+bcf +be^2+afg \cdots(\rm{iv})\\
 a^2e+be^2+bcf+afg=a^2g+bch+aeh+bg^2\cdots(\rm{v})\\
acg+cgh+cdh=ace+cdf+cef\cdots(\rm{vi})\\
ace+cef+cdf=acg+cdh+cgh\cdots(\rm{vii})\\
abc+bd^2+adf+bch=be^2+df^2+ef^2+efh\cdots(\rm{viii})\\
abc+bd^2+bcf+adh= bg^2+dh^2+fgh+gh^2\cdots(\rm{ix})\\
abg+bgh+bdh=abe+bdf+bef\cdots(\rm{xi})\\
abe+bef+bdf=abg+bdh+bgh\cdots(\rm{xii})\\
bcg+dfg+d^2h+ch^2=bce+d^2f+cf^2+deh\cdots(\rm{xiii})\\
bce+deh+d^2f+cf^2=bcg+d^2h+dfg+ch^2\cdots(\rm{xiv}).    
\end{array} \right. 
\end{equation}
\vspace{3mm}

We see easily that (iv)$\Leftrightarrow$(v), (vi)$\Leftrightarrow$(vii),  (xi)$\Leftrightarrow$(xii) and (xiii)$\Leftrightarrow$(xiv), and hence (8) can be rewritten as
\begin{equation}
\left\{\begin{array}{@{\,}lll}   
a^2c+bcd+ade+bcg=ae^2+cf^2+e^2f+efg\\
a^2c+bcd+bce+adg=ag^2+ch^2+egh+g^2h\\
a^2g+bg^2+bch+aeh=a^2e+bcf+be^2+afg\\
c(ag+gh+dh)=c(ae+df+ef)\\
abc+bd^2+adf+bch=be^2+df^2+ef^2+efh\\
abc+bd^2+bcf+adh= bg^2+dh^2+fgh+gh^2\\
b(ag+gh+dh)=b(ae+df+ef)\\
bcg+dfg+d^2h+ch^2=bce+d^2f+cf^2+deh.   
\end{array} \right. 
\end{equation}
\vspace{3mm}
Therefore (9) can be rewritten as
\begin{equation}
\left\{\begin{array}{@{\,}lll}   
a_1^2a_2+b_1a_2b_2+a_1b_2a_3+b_1a_2a_4=a_1a_3^2+a_2b_3^2+a_3^2b_3+a_3b_3a_4\\
a_1^2a_2+b_1a_2b_2+b_1a_2a_3+a_1b_2a_4=a_1a_4^2+a_2b_4^2+a_3a_4b_4+a_4^2b_4\\
a_1^2a_4+b_1a_4^2+b_1a_2b_4+a_1a_3b_4=a_1^2a_3+b_1a_2b_3+b_1a_3^2+a_1b_3a_4\\
a_2(a_1a_4+a_4b_4+b_2b_4)=a_2(a_1a_3+b_2b_3+a_3b_3)\\
a_1b_1a_2+b_1b_2^2+a_1b_2b_3+b_1a_2b_4=b_1a_3^2+b_2b_3^2+a_3b_3^2+a_3b_3b_4\\
a_1b_1a_2+b_1b_2^2+b_1a_2b_3+a_1b_2b_4= b_1a_4^2+b_2b_4^2+b_3a_4b_4+a_4b_4^2\\
b_1(a_1a_4+a_4b_4+b_2b_4)=b_1(a_1a_3+b_2b_3+a_3b_3)\\
b_1a_2a_4+b_2b_3a_4+b_2^2b_4+a_2b_4^2=b_1a_2a_3+b_2^2b_3+a_2b_3^2+b_2a_3b_4.  
\end{array} \right. 
\end{equation}
That is, (5) holds for all $x, y\in A$ iff the scalars $a_i, b_i\, \, (1\le i\le 4)$ satisfy (10).  Therefore we have the following.
\begin{pro}
Let $A$ be a 2-dimensional algebra $A$ over $K$ with structure matrix $\small{\begin{pmatrix}a_1&b_1\\a_2&b_2\\a_3&b_3\\a_4&b_4\end{pmatrix}}$.  Then $A$ is endo-commutative iff  the scalars $a_i, b_i\, \, (1\le i\le 4)$ satisfy $\rm{(10)}$.
\end{pro}

\section{Curled algebras}\label{sec:curled}
Recall that a 2-dimensional algebra is curled if the square of any element $x$ is a scalar multiplication of $x$, otherwise it is straight.  Let $A$ be a two-dimensional algebra over a non-trivial field $K$ with linear base $\{e, f\}$.  If $A$ is curled, we may assume that $e^2=\varepsilon e$ and $f^2=\delta f$ by replacing the bases, where $\varepsilon, \delta\in\{0, 1\}$  .  Write $ef=ae+bf$ and $fe=ce+df$, hence the structure matrix of $A$ is 
\[
C(a, b, c, d; \varepsilon, \delta)\equiv\begin{pmatrix}\varepsilon&0\\0&\delta\\a&b\\c&d\end{pmatrix}, 
\]
where $\varepsilon, \delta\in\{0, 1\}$ and $a, b, c,d\in K$.  However the algebra $C(a, b, c, d; \varepsilon, \delta)$ is not necessarily curled.  The following lemma gives a necessary and sufficient condition for $C(a, b, c, d; \varepsilon, \delta)$ to be curled. 

\begin{lem}
The algebra $C(a, b, c, d; \varepsilon, \delta)$ is curled iff $\small{\left\{\begin{array}{@{\,}lll}   
b+d=\varepsilon\\
a+c=\delta
\end{array} \right. }$ holds.
\end{lem}

\begin{proof}
 Let $x=x_1e+x_2f$ be an arbitrary element of $C(a, b, c, d; \varepsilon, \delta)$, where $x_1, x_2\in K$.  Note that
 \begin{equation}
x^2=\{x_1^2\varepsilon+x_1x_2(a+c)\}e+\{x_2^2\delta+x_1x_2(b+d)\}f.
\end{equation}
 Suppose that $C(a, b, c, d; \varepsilon, \delta)$ is curled.  Then there is $\lambda=\lambda(x_1, x_2)\in K$ such that $x^2=\lambda x$.  Therefore we have from (11) that $\small{\left\{\begin{array}{@{\,}lll}   
x_1^2\varepsilon+x_1x_2(a+c)=\lambda x_1\\
x_2^2\delta+x_1x_2(b+d)=\lambda x_2.
\end{array} \right. }$  By taking $x_1=x_2=1$, we obtain $\varepsilon+a+c=\delta+b+d\, (=\lambda(1, 1))$.  Since $K\ncong\mathbb{F}_2$, we can take $k\in K$ with $k\ne0, 1$.  Then, taking $x_1=1, x_2=k$, we obtain  $\small{\left\{\begin{array}{@{\,}lll}   
\varepsilon+k(a+c)=\lambda(1, k)\\
k^2\delta+k(b+d)=\lambda(1, k)k,
\end{array} \right. }$which implies $\varepsilon+k(a+c)=k\delta+b+d$ because $k\ne0$.  Then we obtain  $\small{\left\{\begin{array}{@{\,}lll}   
\varepsilon+a+c=\delta+b+d\\\varepsilon+k(a+c)=k\delta+b+d,
\end{array} \right. }$  which implies $\small{\left\{\begin{array}{@{\,}lll}   
b+d=\varepsilon\\
a+c=\delta
\end{array} \right. }$because $k\ne1$.

Conversely, suppose that  $\small{\left\{\begin{array}{@{\,}lll}   
b+d=\varepsilon\\
a+c=\delta.
\end{array} \right. }$  Define $\lambda=x_1\varepsilon+x_2\delta$.  By (11), we have $x^2=\lambda x_1e+\lambda x_2f=\lambda x$.  Since $x\in C(a, b, c, d; \varepsilon, \delta)$ is arbitrary, it follows that $C(a, b, c, d; \varepsilon, \delta)$ is curled.  
\end{proof}

Now take $a_1=\varepsilon, b_1=0, a_2=0, b_2=\delta, a_3=a, b_3=b, a_4=c$ and $b_4=d$ in Proposition 2.  In this case, (10) can be rewritten as
\begin{equation}
\left\{\begin{array}{@{\,}lll}   
\varepsilon \delta a=\varepsilon a^2+a^2b+abc\\
\varepsilon \delta c=\varepsilon c^2+acd+c^2d\\
\varepsilon c+\varepsilon ad=\varepsilon a+\varepsilon bc\\
\varepsilon \delta b=\delta b^2+ab^2+abd\\
\varepsilon \delta d=\delta d^2+bcd+cd^2\\
\delta bc+\delta d=\delta b+\delta ad.   
\end{array} \right. 
\end{equation}
Then by Proposition 2, we have the following.

\begin{cor}
The algebra $C(a, b, c, d; \varepsilon, \delta)$ is endo-commutative iff  the scalars $a, b, c, d, \varepsilon$ and $\delta$ satisfy $\rm{(12)}$.
\end{cor}

Moreover, by Lemma 3 and Corollary 2, we have the following.

\begin{lem}
The algebra $C(a, b, c, d; \varepsilon, \delta)$ is endo-commutative and curled iff
\begin{equation}
\left\{\begin{array}{@{\,}lll}   
b+d=\varepsilon \\
a+c=\delta \\
\varepsilon \delta a=\varepsilon a^2+a^2b+abc \\
\varepsilon \delta c=\varepsilon c^2+acd+c^2d \\
\varepsilon c+\varepsilon ad=\varepsilon a+\varepsilon bc \\
\varepsilon \delta b=\delta b^2+ab^2+abd\\
\varepsilon \delta d=\delta d^2+bcd+cd^2 \\
\delta bc+\delta d=\delta b+\delta ad\\   
\end{array} \right. 
\end{equation}
holds.
\end{lem}

We define the following four families of two-dimensional algebras over $K$:
\[
\left\{\begin{array}{@{\,}lll} 
\mathcal{ECC}_{00}:=\{C(a, b, -a, -b; 0, 0) : a, b\in K\}\\
\mathcal{ECC}_{10}:=\{C(0, b, 0, 1-b; 1, 0) : b\in K\}\\
\mathcal{ECC}_{01}:=\{C(a, 0, 1-a, 0; 0, 1) : a\in K\}\\
\mathcal{ECC}_{11}:=\{C(a, 1-a, 1-a, a; 1, 1) : a\in K\}.
\end{array} \right. 
\]
\vspace{1mm}

\noindent
Then we have the following.

\begin{pro}
All algebras in $\mathcal {ECC}_{00}\cup\mathcal{ECC}_{10}\cup\mathcal {ECC}_{10}\cup\mathcal {ECC}_{11}$ are endo-commutative and curled. Conversely, an arbitrary endo-commutative curled algebra of dimension 2 over $K$ is isomorphic to either one of algebras in $\mathcal {ECC}_{00}\cup\mathcal{ECC}_{10}\cup\mathcal {ECC}_{10}\cup\mathcal {ECC}_{11}$. 
\end{pro}

\begin{proof}
If $\varepsilon=\delta=0$, then we see easily that (13)$\Leftrightarrow b+d=a+c=0$, hence we see from Lemma 4 that all algebras in $\mathcal {ECC}_{00}$ are endo-commutative and curled.  If $\varepsilon=1$ and $\delta=0$, then we see easily that (13) $\Leftrightarrow\left\{\begin{array}{@{\,}lll}   
b+d=1\\a=c=0,\end{array}\right.$ hence we see from Lemma 4 that all algebras in $\mathcal {ECC}_{10}$ are endo-commutative and curled.  If $\varepsilon=0$ and $\delta=1$, then we see easily that (13) $\Leftrightarrow\left\{\begin{array}{@{\,}lll}a+c=1\\b=d=0,\end{array}\right.$ hence we see from Lemma 4 that all algebras in $\mathcal {ECC}_{01}$ are endo-commutative and curled.   If $\varepsilon=1$ and $\delta=1$, then we see easily that (13) $\Leftrightarrow b+d=a+c=a+b=1$, hence we see from Lemma 4 that all algebras in $\mathcal {ECC}_{11}$ are endo-commutative and curled.

Conversely let $A$ be an arbitrary endo-commutative curled algebra of dimension 2 over $K$.  Then there are $a_0, b_0, c_0, d_0\in K$ and $\varepsilon_0, \delta_0\in\{0, 1\}$ such that $A$ is isomorphic to the algebra $C(a_0, b_0, c_0, d_0; \varepsilon_0, \delta_0)$.  Then, by Lemma 4, the six scalars $a_0, b_0, c_0, d_0, \varepsilon_0$ and $\delta_0$ must satisfy (13).  Therefore $C(a_0, b_0, c_0, d_0; \varepsilon_0, \delta_0)$ must be in  $\mathcal {EC}_{00}\cup\mathcal{EC}_{10}\cup\mathcal {EC}_{10}\cup\mathcal {EC}_{11}$, as observed in the above argument.
\end{proof}
\vspace{2mm}

\section{Curled algebras of dimension 2: unital, commutative and associative cases}\label{sec:3special-curled}

In this section, we determine unital, commutative, or associative curled algebras of dimension 2 over a non-trivial field $K$. 
\vspace{2mm}

(I) Unital case.  Note that a necessary and sufficient condition for $C(a, b, c, d; \varepsilon, \delta)$ to be unital is 
 \[
\exists u\in C(a, b, c, d; \varepsilon, \delta) : ue = eu = e\, \, {\rm{and}}\, \,  uf = fu = f.
\]
Then we see easily that $C(a, b, c, d; \varepsilon, \delta)$ is unital iff there are $\alpha, \beta\in K$ such that
\[
(\sharp)\, \,  \left\{\begin{array}{@{\,}lll}\alpha\varepsilon+\beta c=1\cdots{\rm{(i)}}\\\beta d=0\cdots{\rm{(ii)}}\\\alpha\varepsilon+\beta a=1\cdots{\rm{(iii)}}\\\beta b=0\cdots{\rm{(iv)}}\\\alpha a=0\cdots{\rm{(v)}}\\\alpha b+\beta\delta=1\cdots{\rm{(vi)}}\\\alpha c=0\cdots{\rm{(vii)}}\\\alpha d+\beta \delta=1.\cdots{\rm{(viii)}}\end{array}\right.
\]
If $\beta\ne0$ and $\alpha\ne0$, then $a=b=c=d=0$ from (ii), (iv), (v) and (vii), hence we can easily see that $(\sharp)\Leftrightarrow\left\{\begin{array}{@{\,}lll}a=b=c=d=0\\\delta=\varepsilon=\alpha=\beta=1.\end{array}\right.$ If $\beta\ne0$ and $\alpha=0$, then $d=b=0$ from (ii), (iv), and so $\delta=\beta=1$ from (vi), hence we can easily see that  $(\sharp)\Leftrightarrow\left\{\begin{array}{@{\,}lll}a=c=\delta=\beta=1\\b=d=0\\\varepsilon=0, 1.\end{array}\right.$  If $\beta=0$, then $\alpha=\varepsilon=1$ from (i), hence we can easily see that $(\sharp)\Leftrightarrow\left\{\begin{array}{@{\,}lll}\alpha=\varepsilon=b=d=1\\a=c=0\\\delta=0, 1.\end{array}\right.$  In other words, in order for $C(a, b, c, d; \varepsilon, \delta)$ to be unital, it is necessary and sufficient that any one of $\left\{\begin{array}{@{\,}lll}a=b=c=d=0\\\delta=\varepsilon=1,\end{array}\right.$ $\left\{\begin{array}{@{\,}lll}a=c=\delta=1\\b=d=0\\\varepsilon=0, 1\end{array}\right.$ and $\left\{\begin{array}{@{\,}lll}\varepsilon=b=d=1\\a=c=0\\\delta=0, 1\end{array}\right.$ holds.  Therefore we have the following.

\begin{lem}
{\rm{(i)}} $C(a, b, c, d; \varepsilon, \delta)$ is unital iff $(a, b, c, d; \varepsilon, \delta)$ is equal to either one of 
\[
(0, 1, 0, 1; 1, 0), (1, 0, 1, 0; 0, 1), (0, 0, 0, 0; 1, 1), (1, 0, 1, 0; 1, 1), (0, 1, 0, 1; 1, 1).
\]

{\rm{(ii)}} If $C(a, b, c, d; \varepsilon, \delta)$ is unital, then it is necessarily endo-commutative.
\vspace{2mm}

{\rm{(iii)}} If $C(a, b, c, d; \varepsilon, \delta)$ is curled, then it is non-unital.

\end{lem}

\begin{proof}
The assertion (i) follows from the above argument.  The assertion (ii) follows from (i) and Corollary 2.  The assertion (iii) follows from (i) and Lemma 3.
\end{proof}
\vspace{2mm}

(II) Commutative case.  Let $x=x_1e+x_2f$ and $y=y_1e+y_2f$ be any elements in $C(a, b, c, d; \varepsilon, \delta)$, where $x_1, x_2, y_1, y_2\in K$.  Then we have
\[
\left\{\begin{array}{@{\,}lll}
xy=(x_1y_1\varepsilon+x_1y_2a+x_2y_1c)e+(x_2y_2\delta+x_1y_2b+x_2y_1d)f\\
yx=(y_1x_1\varepsilon+y_1x_2a+y_2x_1c)e+(y_2x_2\delta+y_1x_2b+y_2x_1d)f.
\end{array}\right.
\]
Therefore we see that the algebra $C(a, b, c, d; \varepsilon, \delta)$ is commutative iff
\[
\left\{\begin{array}{@{\,}lll}
x_1y_1\varepsilon+x_1y_2a+x_2y_1c=y_1x_1\varepsilon+y_1x_2a+y_2x_1c\\
x_2y_2\delta+x_1y_2b+x_2y_1d=y_2x_2\delta+y_1x_2b+y_2x_1d
\end{array}\right.
\]
holds for all $x_1, x_2, y_1, y_2\in K$.  We can easily see that the above conditions hold for all $x_1, x_2, y_1, y_2\in K$ iff $a=c$ and $b=d$.  Therefore we have the following.
\begin{lem}
{\rm{(i)}} $C(a, b, c, d; \varepsilon, \delta)$ is commutative iff $a=c$ and $b=d$ hold.
\vspace{2mm}

{\rm{(ii)}} If  ${\rm{char}}\, K\ne2$, then $C(a, b, c, d; \varepsilon, \delta)$ is curled and commutative iff $(a, b, c, d, \varepsilon, \delta)$ is equal to either one of 
\[
(0, 0, 0, 0; 0, 0), (0, 1/2, 0, 1/2; 1, 0), (1/2, 0, 1/2, 0; 0, 1), (1/2, 1/2, 1/2, 1/2; 1, 1).
\] 
If ${\rm{char}}\, K=2$, then $C(a, b, c, d; \varepsilon, \delta)$ is curled and commutative iff $a=c, b=d$ and $\varepsilon=\delta=0$ hold.
\vspace{2mm}

{\rm{(iii)}}  If $C(a, b, c, d; \varepsilon, \delta)$ is curled and commutative, then it is necessarily endo-commutative.
\end{lem}

\begin{proof}
The assertion (i) follows from the above argument.  The assertions (ii) follow from (i) and Lemma 3. The assertion (iii) follows from (ii) and Corollary 2.
\end{proof}
\vspace{1mm}

(III) Associative case.  Put $\left\{\begin{array}{@{\,}lll} x=x_1e+x_2f\\y=y_1e+y_2f\\z=z_1e+z_2f,
\end{array} \right.$and then we have $xy=X_1e+X_2f$, where  
\[
X_1=x_1y_1\varepsilon+x_1y_2a+x_2y_1c\, \, {\rm{and}}\, \,X_2=x_2y_2\delta+x_1y_2b+x_2y_1d.
\]
Similarly we have $yz=Y_1e+Y_2f$, where
\[
Y_1=y_1z_1\varepsilon+y_1z_2a+y_2z_1c\, \, {\rm{and}}\, \, Y_2=y_2z_2\delta+y_1z_2b+y_2z_1d.
\]
Then we have
\[
(xy)z=(X_1z_1\varepsilon+X_1z_2a+X_2z_1c)e+(X_1z_2b+X_2z_1d+X_2z_2\delta)f
\]
and
\[
x(yz)=(x_1Y_1\varepsilon+x_1Y_2a+x_2Y_1c)e+(x_1Y_2b+x_2Y_1d+x_2Y_2\delta)f.
\]
Then $A$ is associative iff
\[
\left\{\begin{array}{@{\,}lll} 
X_1z_1\varepsilon+X_1z_2a+X_2z_1c=x_1Y_1\varepsilon+x_1Y_2a+x_2Y_1c\cdots(\sharp_1)\\
X_1z_2b+X_2z_1d+X_2z_2\delta=x_1Y_2b+x_2Y_1d+x_2Y_2\delta\cdots(\sharp_2)
\end{array} \right. 
\]
holds for all $x_i, y_i, z_i\in\mathbb Z_2\, \, (1\le i\le2)$.  Put
\[
\left\{\begin{array}{@{\,}lll}
    Z_1=x_1y_1z_1\\
    Z_2=x_1y_1z_2\\
    Z_3=x_1y_2z_1\\
    Z_4=x_1y_2z_2\\
    Z_5=x_2y_1z_1\\
    Z_6=x_2y_1z_2\\
    Z_7=x_2y_2z_1\\
    Z_8=x_2y_2z_2.
\end{array}\right. 
\]
\begin{lem}
The eight polynomials $Z_i\, \, (1\le i\le8)$ on $K$ are linearly independent.
\end{lem}

\begin{proof}
Straightforward.
\end{proof}

Note that $(\sharp_1)$ is rewritten by
\begin{align*}
&Z_1\varepsilon+Z_2\varepsilon a+Z_3(a\varepsilon+bc)+Z_4a^2+Z_5(c\varepsilon+dc)+Z_6ca+Z_7\delta c\\
&=Z_1\varepsilon+Z_2(a\varepsilon+ba)+Z_3(c\varepsilon+da)+Z_4\delta a+Z_5\varepsilon c+Z_6ac+Z_7c^2.
\end{align*}
Then we see easily from Lemma 7 that $(\sharp_1)$ holds for all $x_i, y_i, z_i\in K\, \, (1\le i\le2)$ iff 
\[
\left\{\begin{array}{@{\,}lll}
ab=0\\a\varepsilon+bc=c\varepsilon+ad\\a^2=\delta a\\dc=0\\\delta c=c^2.
\end{array}\right. 
\]
Moreover, note that $(\sharp_2)$ is rewritten by
\begin{align*}
&Z_2\varepsilon b+Z_3bd+Z_4(ab+b\delta)+Z_5d^2+Z_6(cb+d\delta)+Z_7\delta d+Z_8\delta \\
&=Z_2b^2+Z_3bd+Z_4\delta b+Z_5\varepsilon d+Z_6(ad+b\delta)+Z_7(cd+d\delta)+Z_8\delta.
\end{align*}
Then we see easily from Lemma 7 that $(\sharp_2)$ holds for all $x_i, y_i, z_i\in K\, \, (1\le i\le2)$ iff 
\[
\left\{\begin{array}{@{\,}lll}
\varepsilon b=b^2\\ab=0\\d^2=\varepsilon d\\cb+d\delta=ad+b\delta\\cd=0.
\end{array}\right. 
\]
Therefore, $A$ is associative iff
\[
(\sharp)\left\{\begin{array}{@{\,}lll}
ab=0\\a\varepsilon+bc=c\varepsilon+ad\\a^2=\delta a\\dc=0\\\delta c=c^2\\\varepsilon b=b^2\\d^2=\varepsilon d\\cb+d\delta=ad+b\delta.
\end{array}\right. 
\]

The case where $\delta=0$.  In this case, we see that $(\sharp)$ $\Leftrightarrow$ $\left\{\begin{array}{@{\,}lll}
a=0\\c=0\\\varepsilon b=b^2\\\varepsilon d=d^2.
\end{array}\right.$Then the solutions $(a, b, c, d; \varepsilon, \delta)$ of $(\sharp)$ are
\[
(0, 0, 0, 0; 0, 0), (0, 0, 0, 0; 1, 0), (0, 0, 0, 1; 1, 0), (0, 1, 0, 0; 1, 0), (0, 1, 0, 1; 1, 0).
\]

The case where $\delta=1$.  In this case, we see that 
\[
(\sharp)\Leftrightarrow
\left\{\begin{array}{@{\,}lll}
ab=0\\a\varepsilon+bc=c\varepsilon+ad\\a^2=a\\dc=0\\c^2=c\\\varepsilon b=b\\d^2=\varepsilon d\\cb+d=ad+b.
\end{array}\right. 
\]
Then the solutions $(a, b, c, d; \varepsilon, \delta)$ of $(\sharp)$ are
\begin{align*}
&(0, 0, 0, 0; 0, 1), (0, 0, 0, 0; 1, 1),  (0, 1, 0, 1; 1, 1)\cdots{\rm(i)},\\
&(0, 0, 1, 0; 0, 1), (0, 1, 1, 0; 1, 1)\cdots{\rm(ii)},\\
&(1, 0, 0, 0; 0, 1), (1, 0, 0, 1; 1, 1)\cdots{\rm(iii)},\\ 
&(1, 0, 1, 0; 0, 1), (1, 0, 1, 0; 1, 1)\cdots{\rm(iv)}.
\end{align*}
Then we have the following.
\begin{lem}
{\rm{(i)}} $C(a, b, c, d; \varepsilon, \delta)$ is associative iff $(a, b, c, d; \varepsilon, \delta)$ is equal to either one of 

$(0, 0, 0, 0; 0, 0)$, 

$(0, 0, 0, 0; 1, 0), (0, 0, 0, 1; 1, 0), (0, 1, 0, 0; 1, 0), (0, 1, 0, 1; 1, 0)$, 

$(0, 0, 0, 0; 0, 1), (0, 0, 1, 0; 0, 1), (1, 0, 0, 0; 0, 1), (1, 0, 1, 0; 0, 1)$,

$(0, 0, 0, 0; 1, 1), (0, 1, 0, 1; 1, 1), (0, 1, 1, 0; 1, 1), (1, 0, 0, 1; 1, 1), (1, 0, 1, 0; 1, 1)$
\vspace{2mm}

{\rm{(ii)}} If $C(a, b, c, d; \varepsilon, \delta)$ is associative, then it is necessarily endo-commutative.
\vspace{2mm}

{\rm{(iii)}}  If $C(a, b, c, d; \varepsilon, \delta)$ is curled and associative iff $(a, b, c, d; \varepsilon, \delta)$ is equal to either one of 

$(0, 0, 0, 0; 0, 0)$,

$(0, 0, 0, 1; 1, 0), (0, 1, 0, 0; 1, 0)$,

$(0, 0, 1, 0; 0, 1), (1, 0, 0, 0; 0, 1)$, 

$(0, 1, 1, 0; 1, 1), (1, 0, 0, 1; 1, 1)$.
\end{lem}

\begin{proof}
The assertion (i) follows from the above argument.  The assertion (ii) follows from (i) and Corollary 2.  The assertion (iii) follows from (i) and Lemma 3.
\end{proof}
\vspace{3mm}

The following result immediately follows from Lemmas 5, 6 and 8.

\begin{pro}
Suppose that $A$ is a curled algebra of dimension 2 over $K$.  Then
\vspace{2mm}

${\rm{(i)}}$  $A$ is non-unital.
\vspace{2mm}

${\rm{(ii)}}$ When ${\rm{char}}\, K\ne2$,  $A$ is commutative iff  it is isomorphic to either one of the following four algebras:

$C(0, 0, 0, 0; 0, 0)$, 

$C(0, 1/2, 0, 1/2; 1, 0)$,

 $C(1/2, 0, 1/2, 0; 0, 1)$, 
 
  $C(1/2, 1/2, 1/2, 1/2; 1, 1)$.
\vspace{1mm}
  
\noindent
When ${\rm{char}}\, K=2$, $A$ is commutative iff there are $a, b\in K$ such that $A\cong C(a, b, a, b; 0, 0)$.
\vspace{0.1mm}

${\rm{(iii)}}$  $A$  is associative iff  it is isomorphic to either one of the following seven algebras: 

$C(0, 0, 0, 0; 0, 0)$,

$C(0, 0, 0, 1; 1, 0), C(0, 1, 0, 0; 1, 0)$,

$C(0, 0, 1, 0; 0, 1), C(1, 0, 0, 0; 0, 1)$, 

$C(0, 1, 1, 0; 1, 1), C(1, 0, 0, 1; 1, 1)$.

\end{pro}

\begin{cor}
Suppose that $A$ is a curled algebra of dimension 2 over $K$.  If $A$ is either commutative or associative, then $A$ is necessarily endo-commutative.
\end{cor}

\begin{proof}
This follows directly from Lemma 6 (ii) and Lemma 8 (ii).
\end{proof}
\vspace{2mm}

\section{Classification of endo-commutative curled algebras of dimension 2}\label{sec:class-ec-curled}
In this section, we classify endo-commutative curled algebras of dimension 2 over a non-trivial field $K$ by investigating isomorphism of each pair of algebras appearing in $\mathcal{ECC}_{00}\cap \mathcal{ECC}_{10}\cap\mathcal{ECC}_{01}\cap \mathcal{ECC}_{11}$. 
\vspace{2mm}

(I) The classification of algebras in $\mathcal{ECC}_{00}$.  Let $A$ and $A'$ be endo-commutative curled algebras defined by $C(a, b, -a, -b; 0, 0)$ and $C(a', b', -a', -b'; 0, 0)$, respectively.  By Proposition 1, we see that $A\cong A'$ iff there is a nonsingular matrix $X=\begin{pmatrix}x&y\\z&w\end{pmatrix}$ with $\widetilde{X}A'=AX$, which is rewritten as
\[
\begin{pmatrix}x^2&y^2&xy&xy\\z^2&w^2&zw&zw\\xz&yw&xw&yz\\xz&yw&yz&xw\end{pmatrix}\begin{pmatrix}0&0\\0&0\\a'&b'\\-a'&-b'\end{pmatrix}=\begin{pmatrix}0&0\\0&0\\a&b\\-a&-b\end{pmatrix}\begin{pmatrix}x&y\\z&w\end{pmatrix},
\]
that is,
\[
\begin{pmatrix}0&0\\0&0\\a'xw-a'yz&b'xw-b'yz\\
a'yz-a'xw&b'yz-b'xw\end{pmatrix}=\begin{pmatrix}0&0\\0&0\\ax+bz&ay+bw\\-ax-bz&-ay-bw\end{pmatrix},
\]
which is rewritten as \small{$\left\{\begin{array}{@{\,}lll}   
a'(xw-yz)=ax+bz\\b'(xw-yz)=ay+bw
\end{array} \right.$}, that is, $X^t\begin{pmatrix}a\\b\end{pmatrix}=|X|\begin{pmatrix}a'\\b'\end{pmatrix}$.  Note that if $(a, b)\ne(0, 0)$ and $(a', b')\ne(0, 0)$, then we can find a nonsingular $2\times2$ matrix $X$ over $K$ such that $X^t\begin{pmatrix}a\\b\end{pmatrix}=|X|\begin{pmatrix}a'\\b'\end{pmatrix}$, hence $A\cong A'$.  Thus we see that $\mathcal{ECC}_{00}$ is classified into two algebras $C_0$ and $C_1$ defined by $C(0, 0, 0, 0; 0, 0)$ and $C(1, 0, -1, 0; 0, 0)$, respectively up to isomorphism:
\[
C_0=\begin{pmatrix}0&0\\0&0\\0&0\\0&0\end{pmatrix}=\begin{pmatrix}0&0\\0&0\end{pmatrix}\, \, {\rm{and}}\, \, C_1=\begin{pmatrix}0&0\\0&0\\1&0\\-1&0\end{pmatrix}=\begin{pmatrix}0&e\\-e&0\end{pmatrix}.
\]
\vspace{0mm}

II) The classification of algebras in $\mathcal{ECC}_{10}$. Let $A$ and $A'$ be endo-commutative curled algebras defined by $C(0, b, 0, 1-b; 1, 0)$ and $C(0, b', 0, 1-b'; 1, 0)$, respectively.  By Proposition 1, we see that $A\cong A'$ iff there is a nonsingular matrix $X=\begin{pmatrix}x&y\\z&w\end{pmatrix}$ with $\widetilde{X}A'=AX$, which is rewritten as
\[
\begin{pmatrix}x^2&y^2&xy&xy\\z^2&w^2&zw&zw\\xz&yw&xw&yz\\xz&yw&yz&xw\end{pmatrix}\begin{pmatrix}1&0\\0&0\\0&b'\\0&1-b'\end{pmatrix}=\begin{pmatrix}1&0\\0&0\\0&b\\0&1-b\end{pmatrix}\begin{pmatrix}x&y\\z&w\end{pmatrix},
\]
that is,
\[
\begin{pmatrix}x^2&xy\\z^2&zw\\xz&b'xw+(1-b')yz\\
xz&b'yz+(1-b')xw\end{pmatrix}=\begin{pmatrix}x&y\\0&0\\bz&bw\\(1-b)z&(1-b)w\end{pmatrix},
\]
which is rewritten as
\[
\left\{\begin{array}{@{\,}lll}   
x^2=x\\xy=y\\z=0\\b'xw=bw\\(1-b')xw=(1-b)w,
\end{array} \right.
\]
which is equivalent to $[x=1, z=0, b'=b, w\ne0]$ because $|X|\ne0$.  Therefore we see that $A\cong A'$ iff $b'=b$.  Thus we see that $\mathcal{ECC}_{10}$ is classified into the algebras $\{C_2(a)\}_{a\in K}$ defined by $C(0, a, 0, 1-a; 1, 0)\, \, (a\in K)$ up to isomorphism:
\[
C_2(a)=\begin{pmatrix}1&0\\0&0\\0&a\\0&1-a\end{pmatrix}=\begin{pmatrix}e&af\\(1-a)f&0\end{pmatrix}\, \, (a\in K).
\]
\vspace{3mm}

(III) The classification of algebras in $\mathcal{ECC}_{01}$. Let $A$ and $A'$ be endo-commutative curled algebras defined by $C(a, 0, 1-a, 0; 0, 1)$ and $C(a', 0, 1-a', 0; 0, 1)$, respectively.  By Proposition 1, we see that $A\cong A'$ iff there is a nonsingular matrix $X=\small{\begin{pmatrix}x&y\\z&w\end{pmatrix}}$ with $\widetilde{X}A'=AX$, which is rewritten as
\[
\begin{pmatrix}x^2&y^2&xy&xy\\z^2&w^2&zw&zw\\xz&yw&xw&yz\\xz&yw&yz&xw\end{pmatrix}\begin{pmatrix}0&0\\0&1\\a'&0\\1-a'&0\end{pmatrix}=\begin{pmatrix}0&0\\0&1\\a&0\\1-a&0\end{pmatrix}\begin{pmatrix}x&y\\z&w\end{pmatrix},
\]
that is,
\[
\begin{pmatrix}xy&y^2\\zw&w^2\\a'xw+(1-a')yz&yw\\
a'yz+(1-a')xw&yw\end{pmatrix}=\begin{pmatrix}0&0\\z&w\\ax&ay\\(1-a)x&(1-a)y\end{pmatrix},
\]
which is rewritten as
\[
\left\{\begin{array}{@{\,}lll}   
y=0\\zw=z\\w^2=w\\a'xw=ax\\(1-a')xw=(1-a)x,
\end{array} \right.
\]
which is equivalent to $[x\ne0, y=0, w=1, a'=a]$ because $|X|\ne0$.  Therefore we see that $A\cong A'$ iff $a'=a$.  Thus we see that $\mathcal{ECC}_{10}$ is classified into the algebras $\{C_3(a)\}_{a\in K}$ defined by  $C(a, 0, 1-a, 0; 0, 1)\, \,( a\in K)$ up to isomorphism:
\[
C_3(a)=\begin{pmatrix}0&0\\0&1\\a&0\\1-a&0\end{pmatrix}=\begin{pmatrix}0&ae\\(1-a)e&f\end{pmatrix}\, \, (a\in K).
\]
\vspace{0mm}

(IV) The classification of algebras in $\mathcal{ECC}_{11}$. Let $A$ and $A'$ be endo-commutative curled algebras defined by $C(a, 1-a, 1-a, a; 1, 1)$ and  $C(a', 1-a', 1-a', a'; 1, 1)$, respectively.  By Proposition 1, we see that $A\cong A'$ iff there is a nonsingular matrix $X=\small{\begin{pmatrix}x&y\\z&w\end{pmatrix}}$ with $\widetilde{X}A'=AX$, which is rewritten as
\[
\begin{pmatrix}x^2&y^2&xy&xy\\z^2&w^2&zw&zw\\xz&yw&xw&yz\\xz&yw&yz&xw\end{pmatrix}\begin{pmatrix}1&0\\0&1\\a'&1-a'\\1-a'&a'\end{pmatrix}=\begin{pmatrix}1&0\\0&1\\a&1-a\\1-a&a\end{pmatrix}\begin{pmatrix}x&y\\z&w\end{pmatrix},
\]
that is,
\[
\begin{pmatrix}x^2+xy&y^2+xy\\z^2+zw&w^2+zw\\xz+a'xw+(1-a')yz&yw+(1-a')xw+a'yz\\
xz+a'yz+(1-a')xw&yw+(1-a')yz+a'xw\end{pmatrix}=\begin{pmatrix}x&y\\z&w\\ax+(1-a)z&ay+(1-a)w\\(1-a)x+az&(1-a)y+aw\end{pmatrix},
\]
which is rewritten as
\begin{equation}
\left\{\begin{array}{@{\,}lll}   
x^2+xy=x\cdots\cdots(14-1)\\
y^2+xy=y\cdots\cdots(14-2)\\
z^2+zw=z\cdots\cdots(14-3)\\
w^2+zw=w\cdots\cdots(14-4)\\
xz+a'xw+(1-a')yz=ax+(1-a)z\cdots\cdots(14-5)\\
yw+(1-a')xw+a'yz=ay+(1-a)w\cdots\cdots(14-6)\\
xz+a'yz+(1-a')xw=(1-a)x+az\cdots\cdots(14-7)\\
yw+(1-a')yz+a'xw=(1-a)y+aw\cdots\cdots(14-8).
\end{array} \right.
\end{equation}

Suppose that $|X|=xw-yz\ne0$.

(IV-1) The case where $x=0$. Suppose (14) holds. By (14-2), we have $y(y-1)=0$.  If $y=0$, then $|X|=0$, a contradiction.  Then we have $y=1$.  By (14-7), we have $a'z=az$.  But since $|X|=-z\ne0$, we have $a'=a$. 

(IV-2) The case where $x\ne0$ and $z=0$. Suppose (14) holds.  By (14-1), we have $x+y=1$.  By (14-5), we have $a'w=a$.  By (14-7), we have $(1-a')w=1-a$, hence $w=1$ because $a'w=a$.  Thus we see that $a'=a$

(IV-3) The case where $x\ne0$ and $z\ne0$. Suppose (14) holds.  By (14-1), we have $x+y=1$.  By (14-3), we have $z+w=1$.  By (14-5), we have $a'(x-z)=a(x-z)$ because $y=1-x, w=1-z$.  If $x=z$, then $|X|=x(1-z)-(1-x)z=0$, a contradiction. Hence $x\ne z$, so we have $a'=a$. 

Considering the above three cases, we see that $A\cong A'\Leftrightarrow a'=a$, and hence $\mathcal{ECC}_{11}$ is classified into algebras $\{C_4(a)\}_{a\in K}$ defined by $C(a, 1-a, 1-a, a; 1, 1)\, \, (a\in K)$ up to isomorphism:
\[
C_4(a)=\begin{pmatrix}1&0\\0&1\\a&1-a\\1-a&a\end{pmatrix}=\begin{pmatrix}e&ae+(1-a)f\\(1-a)e+af&f\end{pmatrix}\, \, (a\in K).
\]
\vspace{0mm}

(V) Classification between $\mathcal{ECC}_{00}$ and $\mathcal{ECC}_{10}$. Let $A\in\mathcal{ECC}_{00}$ and $A'\in\mathcal{ECC}_{10}$ be arbitrary.   By (I) and (II), we have that ${\rm{rank}}\, A\le1$ and  ${\rm{rank}}\, A'=2$, hence $A\ncong A'$ from Corollary 1.
\vspace{3mm}

(VI) Classification between $\mathcal{ECC}_{00}$ and $\mathcal{ECC}_{01}$. Let $A\in\mathcal{ECC}_{00}$ and $A'\in\mathcal{ECC}_{01}$ be arbitrary.   By (I) and (III), we have that ${\rm{rank}}\, A\le1$ and  ${\rm{rank}}\, A'=2$, hence $A\ncong A'$ from Corollary 1.
\vspace{3mm}

(VII) Classification between $\mathcal{ECC}_{00}$ and $\mathcal{ECC}_{11}$. Let $A\in\mathcal{ECC}_{00}$ and $A'\in\mathcal{ECC}_{11}$ be arbitrary.   By (I) and (IV), we have that ${\rm{rank}}\, A\le1$ and  ${\rm{rank}}\, A'=2$, hence $A\ncong A'$ from Corollary 1.
\vspace{3mm}

(VIII) Classification between $\mathcal{ECC}_{10}$ and $\mathcal{ECC}_{01}$. Let $A\in\mathcal{ECC}_{10}$ and $A'\in\mathcal{ECC}_{01}$ be arbitrary.   By (II) and (III), we have that $A\cong C_2(a)$ and $A'\cong C_3(a')$ for some $a, a'\in K$, where  $C_2(a)$ and $C_3(a')$ are defined by $C(0, a, 0, 1-a; 1, 0)$ and  $C(a', 0, 1-a', 0; 0, 1)$, respectively.

Suppose $C_2(a)\cong C_3(a')$.  By Proposition 1, there is a nonsingular matrix $X=\small{\begin{pmatrix}x&y\\z&w\end{pmatrix}}$ with $\widetilde{X}C(a', 0, 1-a', 0; 0, 1)=C(0, a, 0, 1-a; 1, 0)X$, which is rewritten as
\[
\begin{pmatrix}x^2&y^2&xy&xy\\z^2&w^2&zw&zw\\xz&yw&xw&yz\\xz&yw&yz&xw\end{pmatrix}\begin{pmatrix}0&0\\0&1\\a'&0\\1-a'&0\end{pmatrix}=\begin{pmatrix}1&0\\0&0\\0&a\\0&1-a\end{pmatrix}\begin{pmatrix}x&y\\z&w\end{pmatrix},
\]
that is, 
\[
\begin{pmatrix}xy&y^2\\zw&w^2\\a'xw+(1-a')yz&yw\\a'yz+(1-a')xw&yw\end{pmatrix}=\begin{pmatrix}x&y\\0&0\\az&aw\\(1-a)z&(1-a)w\end{pmatrix},
\]
which is rewritten as
\begin{equation}
\left\{\begin{array}{@{\,}lll}   
xy=x\, \, \, \, \, \, \,\,  \cdots\cdots\cdots\cdots\cdots\cdots\cdots\cdots(15-1)\\
y^2=y\, \, \, \, \, \, \,\,  \cdots\cdots\cdots\cdots\cdots\cdots\cdots\cdots(15-2)\\
zw=0\, \, \, \, \, \, \,\,  \cdots\cdots\cdots\cdots\cdots\cdots\cdots\cdots(15-3)\\
w^2=0\, \, \, \, \, \, \,\,  \cdots\cdots\cdots\cdots\cdots\cdots\cdots\cdots(15-4)\\
a'xw+(1-a')yz=az\, \, \, \, \, \, \cdots\cdots\cdots(14-5)\\
yw=aw\, \, \, \, \cdots\cdots\cdots\cdots\cdots\cdots\cdots\cdots(15-6)\\
a'yz+(1-a')xw=(1-a)z\cdots\cdots(15-7)\\
yw=(1-a)w\, \, \, \, \, \cdots\cdots\cdots\cdots\cdots\cdots(15-8).
\end{array} \right.
\end{equation}
By (15-4), we have $w=0$.  Since $|X|=xw-yz=-yz\ne0$, it follows from (15-2) that $y=1$.  Then we have from (15-7) that $a+a'=1$ because $y=1, w=0$ and $z\ne0$.  Thus we have $C_2(a)\cong C_3(a') \Rightarrow a+a'=1$.  

Conversely, suppose $a+a'=1$.  Define $X=\small{\begin{pmatrix}0&1\\1&0\end{pmatrix}}$, and then 
\[
\widetilde{X}=\begin{pmatrix}0&1&0&0\\1&0&0&0\\0&0&0&1\\0&0&1&0\end{pmatrix},
\]
so we see easily that $\widetilde{X}C(a', 0, 1-a', 0; 0, 1)=C(0, a, 0, 1-a; 1, 0)X$ by the assumption $a+a'=1$.  Therefore we have from Proposition 1 that $A_2(a)\cong A_3(a')$.  Consequently, we see that $C_2(a)\cong C_3(1-a)$ holds for all $a\in K$.
\vspace{3mm}

(IX) Classification between $\mathcal{ECC}_{10}$ and $\mathcal{ECC}_{11}$. Let $A\in\mathcal{ECC}_{10}$ and $A'\in\mathcal{ECC}_{11}$ be arbitrary.   By (II) and (IV), we have that $A\cong C_2(a)$ and $A'\cong C_4(a')$ for some $a, a'\in K$, where  $C_2(a)$ and $C_4(a')$ are defined by $C(0, a, 0, 1-a; 1, 0)$ and  $C(a', 1-a', 1-a', a'; 1, 1)$, respectively.

Suppose $C_2(a)\cong C_4(a')$.  By Proposition 1, there is a nonsingular matrix $X=\small{\begin{pmatrix}x&y\\z&w\end{pmatrix}}$ with $\widetilde{X}C(a', 1-a', 1-a', a'; 1, 1)=C(0, a, 0, 1-a; 1, 0)X$, which is rewritten as
\[
\begin{pmatrix}x^2&y^2&xy&xy\\z^2&w^2&zw&zw\\xz&yw&xw&yz\\xz&yw&yz&xw\end{pmatrix}\begin{pmatrix}1&0\\0&1\\a'&1-a'\\1-a'&a'\end{pmatrix}=\begin{pmatrix}1&0\\0&0\\0&a\\0&1-a\end{pmatrix}\begin{pmatrix}x&y\\z&w\end{pmatrix},
\]
that is, 
\[
\begin{pmatrix}x^2+xy&y^2+xy\\z^2+zw&w^2+zw\\xz+a'xw+(1-a')yz&yw+(1-a')xw+a'yz\\xz+a'yz+(1-a')xw&yw+(1-a')yz+a'xw\end{pmatrix}=\begin{pmatrix}x&y\\0&0\\az&aw\\(1-a)z&(1-a)w\end{pmatrix},
\]
which is rewritten as
\begin{equation}
\left\{\begin{array}{@{\,}lll}   
x^2+xy=x\, \, \,  \, \cdots\cdots\cdots\cdots\cdots\cdots\cdots\cdots(16-1)\\
y^2+xy=y\, \, \, \, \cdots\cdots\cdots\cdots\cdots\cdots\cdots\cdots(16-2)\\
z^2+zw=0\, \, \, \, \,  \cdots\cdots\cdots\cdots\cdots\cdots\cdots\cdots(16-3)\\
w^2+zw=0\, \, \,    \cdots\cdots\cdots\cdots\cdots\cdots\cdots\cdots(16-4)\\
xz+a'xw+(1-a')yz=az\, \, \cdots\cdots\cdots(16-5)\\
yw+(1-a')xw+a'yz=aw\cdots\cdots\cdots(16-6)\\
xz+a'yz+(1-a')xw=(1-a)z\, \, \cdots(16-7)\\
yw+(1-a')yz+a'xw=(1-a)w\cdots(16-8).
\end{array} \right.
\end{equation}

The case where $x=0$. Since $|X|=-yz\ne0$, it follows that $y\ne0$ and $z\ne0$, so $y=1$ and $z+w=0$ by (16-2) and (16-3), respectively.  Then we have $a+a'=1$ by (16-5).

The case where $x\ne0$. By (16-1), we have $x+y=1$.  Since $|X|\ne0$, it follows that either $z\ne0$ or $w\ne0$ holds, hence $z+w=0$ by (16-3) and (16-4).  Substituting $y=1-x$ and $w=-z$ into (16-5), we obtain $z(1-a-a')=0$, hence $a+a'=1$ because $z\ne0$ (otherwise, $z=w=0$, so $|X|=0$, a contradiction).  Thus we have $C_2(a)\cong C_4(a') \Rightarrow a+a'=1$.  

Conversely, suppose $a+a'=1$.  Define $X=\small{\begin{pmatrix}0&1\\1&-1\end{pmatrix}}$, and then 
\[
\widetilde{X}=\begin{pmatrix}0&1&0&0\\1&1&-1&-1\\0&-1&0&1\\0&-1&1&0\end{pmatrix},
\]
so we see easily that  $\widetilde{X}C(a', 1-a', 1-a', a'; 1, 1)=C(0, a, 0, 1-a; 1, 0)X$ by the assumption $a+a'=1$.  Therefore we have from Proposition 1 that $C_2(a)\cong C_4(a')$.  Consequently, we see that $C_2(a)\cong C_4(1-a)$ holds for all $a\in K$.
\vspace{3mm}

(X) Classification between $\mathcal{EC}_{01}$ and $\mathcal{EC}_{11}$. By (VIII) and (IX), we see that $A_3(a)\cong A_4(a)$ holds for all $a\in K$.

Therefore, by (I)$\sim$(X) and Proposition 3, we have the following.
\vspace{3mm}

\begin{thm}
Up to isomorphism, two-dimensional endo-commutative curled algebras over a non-trivial field $K$ with a linear basis $\{e, f\}$ are classified into the algebras
\[
C_0, C_1\, \, {\rm{and}}\, \, \{C_2(a)\}_{a\in K}
\]
defined by multiplication tables
\[
\begin{pmatrix}0&0\\0&0\end{pmatrix}, \begin{pmatrix}0&e\\-e&0\end{pmatrix}\, \, {\rm{and}}\, \, \begin{pmatrix}e&af\\(1-a)f&0\end{pmatrix}\, \, (a\in K),
\]
respectively. 
\end{thm}

\begin{cor}
Suppose that $A$ is an endo-commutative curled algebra of dimension 2 over a non-trivial field $K$.  Then

${\rm{(i)}}$ $A$ is zeropotent iff $A\cong C_0$ or $A\cong C_1$.

${\rm{(ii)}}$ When ${\rm{char}}K=2$, $A$ is commutative iff $A\cong C_0$ or $A\cong C_1$.  When ${\rm{char}}K\ne2$,  $A$ is commutative iff $A\cong C_0$ or $A\cong C_2(1/2)$.

${\rm{(iii)}}$  $A$ is anti-commutative iff $A\cong C_0$ or $A\cong C_1$.

${\rm{(iv)}}$  $A$ is associative iff $A\cong C_0, A\cong C_2(0)$ or $A\cong C_2(1)$.
\end{cor}

\begin{proof}
(i) Note that $C_0$ and $C_1$ are zeropotent, but $C_2(a)\, \, (a\in K)$ is not zropotent.  Then we obtain the desired result from Theorem 1.
\vspace{2mm}

(ii) Suppose  ${\rm{char}}K=2$.  It is obvious that $C_0$ and $C_1$ are commutative.  But since $af\ne(1-a)f$ for all $a\in K$, it follows that  $C_2(a)$ is non-commutative for all $a\in K$.  Then we see from Theorem 1 that $A$ is commutative iff $A\cong C_0$ or $A\cong C_1$.  

Next suppose ${\rm{char}}K\ne2$.  It is obvious that $C_0$ is commutative, but $C_1$ is non-commutative.  Also, we see that $af=(1-a)f$ iff $a=1/2$. Then we see from Theorem 1 that $A$ is commutative iff $A\cong C_0$ or $A\cong C_2(1/2)$.  
\vspace{2mm}

(iii) Note that $C_0$ and $C_1$ are anti-commutative, but $C_2(a)$ is not anti-commutative for all $a\in K)$.  Then we obtain the desired result from Theorem 1.
\vspace{2mm}

(iv) By Proposition 4 (iii), we see that $A$ is associative iff it is isomorphic to any one of $C_0, C_2(0), C_2(1), C_3(0), C_3(1), C_4(0)$ and $C_4(1)$.  However, as observed in (VIII) and (IX), we see that $C_2(a)\cong C_3(1-a)$ and $C_2(a)\cong C_4(1-a)$ for all $a\in K$.  Therefore we obtain the desired result.
\end{proof}

\end{document}